\documentclass[11pt,reqno]{amsart}
\usepackage{pb-diagram} 
\usepackage{amsmath,amssymb}
\usepackage{amsthm}
\usepackage{hyperref}
\usepackage{pb-diagram}

\newtheorem{theorem}{Theorem}

\newtheorem{corollary}[theorem]{Corollary}
\newtheorem{remark}{Remark}
\newtheorem{lemma}[theorem]{Lemma}
\newtheorem{thrmst}{Theorem ([24])}

\newcommand{\eend}{\mathrm{End}\,}

\newcommand{\trace}{\mathrm{Trace\,}}

\newcommand{\C}{\mathbb{C}}

\newcommand{\beqt}{\begin{equation}}  \newcommand{\eeqt}{\end{equation}}
\newcommand{\bal}{\begin{align}}      \newcommand{\eal}{\end{align}}
\newcommand{\ba}{\begin{array}}      \newcommand{\ea}{\end{array}}
\newcommand{\bc}{\begin{center}}     \newcommand{\ec}{\end{center}}
\newcommand{\be}{\begin{enumerate}}  \newcommand{\ee}{\end{enumerate}}
\newcommand{\beq}{\begin{eqnarray}}  \newcommand{\eeq}{\end{eqnarray}}
\newcommand{\beQ}{\begin{eqnarray*}} \newcommand{\eeQ}{\end{eqnarray*}}
\newcommand{\bi}{\begin{itemize}}    \newcommand{\ei}{\end{itemize}}
\newcommand{\bt}{\begin{tabular}}    \newcommand{\et}{\end{tabular}}

\let\<\langle     
\let\>\rangle   

\title{Rigidity of compact Riemannian spin Manifolds with Boundary}
\author{Simon Raulot$^1$}\thanks{$^1$ Supported by the Swiss SNF grant $20$-$118014/1$}
\address
{Institut de Math\'ematiques\\
Universit\'e de Neuch\^atel\\
Rue Emile-Argand 11 \\2007 Neuch\^atel\\ Suisse}
\email{simon.raulot@unine.ch}
\date{\today}


\begin{document}

\maketitle 

\begin{abstract}
In this article, we prove new rigidity results for compact Riemannian spin manifolds with boundary whose scalar curvature is bounded from below by a non-positive constant. In particular, we obtain generalizations of a result of Hang-Wang \cite{hangwang1} based on a conjecture of Schroeder and Strake \cite{schroeder}.
\end{abstract}



\section{Introduction}\label{intro}


The well-known spinorial proof of the positive mass theorem for asymptotically flat manifolds given by Witten \cite{witten} is based on a subtle use of the Weitzenb\"ock type formula for the hypersurface Dirac-type operator. In this setting, asymptotic flatness provides a boundary condition for the metric at infinity. A corollary of the positive mass theorem is that there is no Riemannian metric on $\mathbb{R}^n$ with nonnegative scalar curvature which is euclidean outside a compact set, except the Euclidean one. Using Witten's approach, rigidity results for noncompact manifolds whose metric behaviour is prescribed at infinity were intensively studied (see for example \cite{anderson}, \cite{herzlich} or \cite{minoo}).\\

More recently, rigidity results for compact manifolds with boundary have been proved using generalizations of Witten's positive mass theorem (see \cite{miao} or \cite{shitam} for example). In this setting, the conditions on the metric at infinity are replaced with natural conditions on the metric of the boundary. The latter are tightly related to the notion of quasi-local mass in General Relativity.\\

On the other hand, recent papers by Hijazi, Montiel, Rold\'an and Zhang (see \cite{HMZ1}, \cite{HMR} or \cite{HM}) emphasize the fact that Spin geometry provides an adaptated framework for the study of hypersurfaces. In particular, they show that under intrinsic and extrinsic curvature assumptions on a Riemannian spin manifold with boundary, there is an isomorphism between the restriction to the boundary of parallel spinors and extrinsic Killing spinors. \\

In this paper, we generalize the results of \cite{HM}. We prove that, under suitable assumptions, a solution of the Dirac equation:
\begin{equation}\label{de}
{\bf D }\Phi=\frac{n-1}{2} H_0\Phi,
\end{equation}

can be extend to a parallel spinor field on the whole manifold. Here ${\bf{D}}$ is the extrinsic Dirac operator of the boundary (see Section \ref{Pg}) and $H_0$ is a nonnegative (and non identically zero) function on $\partial M$. Several rigidity results follow by noting that the existence of a spinor field satisfying Equation (\ref{de}) is tightly related to the existence of an isometric immersion of the boundary in a manifold carrying a parallel spinor field. One of these applications is given by a generalization of a theorem in \cite{hangwang1} which improves a conjecture of Schroeder and Strake \cite{schroeder} in the spinorial setting.\\

In the last part of this paper, we study the hyperbolic version of the results obtained in the previous paragraph.\\

{\bf Acknowledgements:} I would like to thank Oussama Hijazi, Emmanuel Humbert and Julien Roth for their remarks and suggestions, as well as their support. I am also grateful to the ``Institut de Math\'ematiques'' of the University of Neuch\^atel for its financial support. Finally, I would like to thank the referees for helpful comments.


\section{Geometric Preliminaries}\label{Pg}


Let $(M^n,g)$ be an $n$-dimensional Riemannian spin manifold. We denote by $\Sigma M$ the bundle of complex spinor fields over $M$ and by $\nabla$ the Riemannian and the spin Levi-Civita connections. The Clifford multiplication, that is the action of the Clifford bundle $\C l(M)$ on the spinor bundle, will be denoted by:
\begin{eqnarray*}
\gamma : \mathbb{C} l (M)\longrightarrow \eend(\Sigma M)
\end{eqnarray*}

\noindent and the natural Hermitian product on $\Sigma M$ compatible with $\nabla$ and $\gamma$ by $\<\; ,\;\>$. The Dirac operator is defined by taking the composition of the Clifford multiplication with the spinorial Levi-Civita connection that is $D=\gamma\circ\nabla$. It is an elliptic differential operator of order one acting on the sections of the spinor bundle and it is locally given by:
\begin{eqnarray*}
D=\sum_{i=1}^{n}\gamma(e_i)\nabla_{e_i},
\end{eqnarray*}

\noindent where $\left\lbrace  e_1,\cdots,e_{n}\right\rbrace$ is a local $g$-orthonormal frame of $TM$.

\noindent Assume now that $M$ has a smooth boundary $\partial M$. Since $\partial M$ is an oriented hypersurface of $M$, its normal bundle is trivial and thus one can define a spin structure on $\partial M$. Then we can build the {\it intrinsic} spinor bundle over $\partial M$ denoted by $\Sigma(\partial M)$ which is naturally endowed with the spin Levi-Civita connection $\nabla^{\partial M}$, a Clifford multiplication $\gamma^{\partial M}$ and also the Dirac operator $D^{\partial M}$, called the {\it intrinsic} Dirac operator of $\partial M$. One can also define (see \cite{Ba2} for example) an {\it extrinsic} spinor bundle over $\partial M$ by putting ${\bf{S}}:=\Sigma M_{|\partial M}$. This bundle is also endowed with the spin Levi-Civita connection $\nabla^{{\bf{S}}}$ and a Clifford multiplication $\gamma^{{\bf{S}}}$ which can be related with these acting on $\Sigma M$ by:
\begin{eqnarray}\label{gauss} 
\nabla_X & = & \nabla^{{\bf{S}}}_X+\frac{1}{2}\gamma^{\bf{S}}\big(A(X)\big)\quad\text{\it{(Spinorial Gauss Formula)}}\\
\gamma^{{\bf{S}}}(X) & = & \gamma(X)\gamma(\nu)\nonumber,
\end{eqnarray}

\noindent for all $X\in\Gamma\big(T(\partial M)\big)$ and where $\nu$ is the inward unit vector field normal to $\partial M$ and $A$ is the (symmetric) Weingarten map given by $A(X)=-\nabla_X\nu$. As for the intrinsic case, one can define a Dirac operator acting on $\bf{S}$ by ${\bf{D}}:=\gamma^{\bf{S}}\circ\nabla^{\bf{S}}$. This operator will be called the {\it extrinsic} Dirac operator of $\partial M$. A straightforward calculation using the spinorial Gauss formula (\ref{gauss}) allows to obtain a relation between this operator and the Dirac operator of $M$, namely:
\beqt\label{Dirac bord} 
{\bf{D}}\psi=\frac{n-1}{2}{H}\psi-\gamma(\nu)D\psi-\nabla_{\nu}\psi.
\eeqt

\noindent for all $\psi\in\Gamma({\bf S})$ and where $H:=\frac{1}{n-1}\trace (A)$ is the mean curvature of $\partial M$ in $M$.

\noindent The hypersurface $\partial M$ has thus two spinor bundles which can be identified in a canonical way. Indeed the extrinsic spinor bundle $({\bf{S}},\gamma^{{\bf{S}}},\nabla^{{\bf{S}}},{\bf{D}})$ is isometric to:
\begin{eqnarray}\label{idimpair}
 \big(\Sigma(\partial M),\gamma^{\partial M},\nabla^{\partial M},D^{\partial M}\big) 
\end{eqnarray}

if $n$ is odd and to:
\begin{eqnarray}\label{idpair}
\big(\Sigma(\partial M)\oplus\Sigma(\partial M),\gamma^{\partial M}\oplus -\gamma^{\partial M},\nabla^{\partial M}\oplus \nabla^{\partial M},D^{\partial M}\oplus -D^{\partial M}\big).
\end {eqnarray}

\noindent if $n$ is even. For more details on these identifications, we refer to \cite{Ba2}, \cite{HMZ1} or \cite{Mo}. On the other hand, using the relation:
\begin{eqnarray}\label{nudirac}
{\bf{D}}\gamma(\nu)=-\gamma(\nu){\bf{D}}
\end{eqnarray}
 
one can easily check that the spectrum of the extrinsic Dirac operator is symmetric with respect to zero. Using these identifications, O. Hijazi and S. Montiel \cite{HM} define the notion of extrinsic Killing spinor, generalizing Killing spinors to the frame of hypersurfaces. More precisely, a spinor field $\varphi\in\Gamma({\bf{S}})$ is an {\it extrinsic} Killing spinor if for all $X\in\Gamma\big(T(\partial M)\big)$, we have:
\begin{eqnarray}\label{killingextrinseque}
\nabla^{{\bf{S}}}_X\varphi=-\alpha\gamma^{{\bf{S}}}(X)\varphi,
\end{eqnarray}
where $\alpha\in\mathbb{R}$. Under some curvature assumptions, the authors prove that the existence of such a spinor field on a hypersurface bounding a compact domain implies that the domain carries a parallel spinor (and hence is Ricci flat) and forces the boundary to be totally umbilical with constant mean curvature. In this article, we study this question for spinors solutions of the Dirac equation.\\

{\bf Convention:} In the following we assume that if $\Sigma^{n-1}$ is a compact manifold which is isometrically immersed in two $n$-dimensional manifolds, then these immersions induce the same spin structure on $\Sigma^{n-1}$.


\section{Domains with positive scalar curvature}\label{euclidien}


In this section, we consider  an $n$-dimensional compact Riemannian spin manifold $(M^n, g)$ with nonnegative scalar curvature $R$. Assume that $M$ has a smooth boundary $\partial M$ which has $p$ connected components $\partial M_j$ with nonnegative mean curvature $H^{(j)}$ for all $1\leq j\leq p$. We first prove that under suitable assumptions on the mean curvature of the boundary, we can extend a spinor field satisfying the Dirac equation (\ref{de}) to a parallel spinor field on $M$. Then we apply this result to obtain several rigidity results for compact Riemannian spin manifolds with boundary. The first result we get is the following:
\begin{theorem}\label{main}
Let $(M^n,g)$ be an $n$-dimensional compact and connected Riemannian spin manifold with smooth boundary $\partial M$. Assume that the scalar curvature of $M$ is nonnegative, that the mean curvature of each connected component $\partial M_j$ of $\partial M$ is nonnegative (and non identically zero). If there exists a smooth spinor field $\Phi\in\Gamma({\bf S}_{j_0})$ satisfying:
\begin{eqnarray}\label{diracequation}
{\bf{D}}\Phi=\frac{n-1}{2} H_0\Phi,
\end{eqnarray}
where $H_0$ is a smooth function on $\partial M_{|j_0}$ such that $0\leq H_0\leq H^{(j_0)}$, then $(M^n, g)$ has a parallel spinor field, the boundary is connected and $H^{(j_0)}=H_0$.
\end{theorem} 

In this theorem, we let ${\bf{S}}_{j_0}:=\Sigma M_{|\partial M_{j_0}}$. The proof of this result relies on the Schr\"odinger-Lichnerowicz formula \cite{Li} which gives a relation between the square of the Dirac operator and the spin Laplacian. More precisely, we have:
\begin{eqnarray}\label{lich}
D^2=\nabla^*\nabla+\frac{R}{4},
\end{eqnarray}
and when integrated over $M$ (see \cite{HMZ1}) yields to:
\begin{eqnarray}\label{reilly}
\int_{M}\big(|\nabla\psi|^2-|D\psi|^2+\frac{1}{4}R|\psi|^2\big)dv =\int_{\partial M} \big(\<{\bf{D}}\psi,\psi\>-\frac{n-1}{2}H|\psi|^2\big)ds
\end{eqnarray}

for all $\psi\in\Gamma(\Sigma M)$. In comparison with the classical Reilly formula on functions (see \cite{Re}), Formula (\ref{reilly}) is called the spinorial Reilly formula. The other key point in the proof of Theorem~\ref{main} is an adaptated choice of boundary conditions for the Dirac operator on $M$. For more details on this subject, we refer to \cite{BW} or \cite{hijazi.montiel.roldan:01}. We consider here the ${\rm MIT}$ condition defined by the pointwise orthogonal projection:
$$\begin{array}{lccl}
{P}^{\pm}: &L^2({\bf{S}}) & \longrightarrow & L^2(V^{\pm})\\
 & \varphi & \longmapsto & \frac{1}{2}(Id\pm i\gamma(\nu))\varphi
\end{array}$$ 

where $V^{\pm}$ is the subbundle of ${\bf S}$ whose fiber is the eigenspace associated with the eigenvalue $\pm 1$ of the involution $i\gamma(\nu)$. One can then check that this map defines an elliptic boundary condition for the Dirac operator $D$ of $M$ and we can prove (see \cite{HMZ2}):

\begin{lemma}\label{existence}
Let $(M^n, g)$ be an $n$-dimensional compact Riemannian spin manifold with smooth boundary $\partial M$, then the map
\begin{eqnarray*}
D:\{\varphi\in H^2_1(\Sigma M)\,:\,
P^\pm\varphi_{|\partial M}=0\}\longrightarrow L^2(\Sigma M)
\end{eqnarray*}
\noindent is invertible.
\end{lemma}

We can now give the proof of our first result. \\

{\it Proof of Theorem~\ref{main}:} Let $\Phi\in\Gamma ({\bf{S}}_{j_0})$ a solution of the Dirac equation (\ref{diracequation}) and we extend this spinor field on $M$ by $\widetilde{\Phi}$ in such a way that it vanishes on the other components of $\partial M$, that is:
\begin{eqnarray}\label{prolonge}
\widetilde{\Phi}_j:=\widetilde{\Phi}_{|\partial M_{j}}=
\left\lbrace
\begin{array}{ll}
\Phi & \text{if } j=j_0\\
0 & \text{if } j\neq j_0.
\end{array}
\right.
\end{eqnarray}

Lemma \ref{existence} ensures the existence of a unique smooth spinor field $\Psi\in\Gamma(\Sigma M)$ satisfying the boundary problem:
\begin{equation}\label{pab}
\left\lbrace
\begin{array}{ll}
D\Psi=0 & \qquad\text{on}\;M\\
P^\pm\Psi_{|\partial M}=P^\pm\widetilde{\Phi}_{|\partial M} & \qquad\text{along}\;\partial M,
\end{array}
\right.
\end{equation}
which by (\ref{prolonge}) gives:
\begin{equation}\label{pab1}
\left\lbrace
\begin{array}{ll}
D\Psi=0 & \qquad\text{on}\;M\\
P^\pm\Psi_{|\partial M_{j_0}}=P^\pm\Phi & \qquad\text{along}\;\partial M_{j_0}\\
P^\pm\Psi_{|\partial M_{j}}=0 & \qquad\text{along}\;\partial M_{j}\;\text{for}\;j\neq j_0.
\end{array}
\right.
\end{equation}

In the sequel, we will denote equally a spinor field on $M$ and its restriction on the boundary. Using the spinorial Reilly formula (\ref{reilly}) and since $R\geq 0$, we get:
\begin{eqnarray}\label{reil}
0\leq\int_{M}\Big(|\nabla\Psi|^2+\frac{1}{4}R|\Psi|^2\Big)dv & = & \sum_{j=1}^p\int_{\partial M_j} \Big(\<{\bf{D}}\Psi,\Psi\> -\frac{n-1}{2}H^{(j)}|\Psi|^2\Big)ds.
\end{eqnarray}

We now prove that the boundary term in the preceding formula is nonpositive. First we observe that since the spinor $\Phi$ satisfies the Dirac equation (\ref{diracequation}) on $\partial M_{j_0}$, we obtain with the help of (\ref{nudirac}) that:
\begin{eqnarray}\label{ppm}
{\bf{D}}(P^\pm\Phi)=\frac{n-1}{2}H_0 P^\mp\Phi.
\end{eqnarray}

On the other hand, for all $\varphi\in\Gamma ({\bf{S}})$, an integration by parts using the symmetry of the Dirac operator ${\bf{D}}$ and the decomposition $\varphi=P^+\varphi+P^-\varphi$ yield:
\begin{eqnarray*}
\int_{\partial M}\<{\bf{D}}\varphi,\varphi\>ds=2\int_{\partial M}{\rm Re}\<{\bf{D}}(P^\pm\varphi),P^\mp\varphi\>ds.
\end{eqnarray*}

Now for $\varphi=\Psi$, we get with (\ref{pab}) and (\ref{ppm}):
\begin{eqnarray}\label{masse}
\int_{\partial M_{j_0}}\<{\bf{D}}\Psi,\Psi\>ds=(n-1)\int_{\partial M_{j_0}}H_0\,{\rm Re}\<P^\mp\Phi,P^\mp\Psi\>ds.
\end{eqnarray}

Moreover, since $|P^\mp\Psi-P^\mp\Phi|^2\geq 0$, we have:
\begin{eqnarray}\label{masse2}
2\,{\rm Re}\<P^\mp\Phi,P^\mp\Psi\>\leq |P^\mp\Psi|^2+|P^\mp\Phi|^2,
\end{eqnarray}

which leads to:
\begin{eqnarray}\label{masse1}
\int_{\partial M_{j_0}}\<{\bf{D}}\Psi,\Psi\>ds\leq\frac{n-1}{2}\int_{\partial M_{j_0}}H_0\big(|P^\mp\Psi|^2+|P^\mp\Phi|^2\big)ds.
\end{eqnarray}

We also remark that the symmetry of ${\bf{D}}$ and (\ref{ppm}) give:
\begin{eqnarray*}
\int_{\partial M_{j_0}}H_0|P^\pm\Phi|^2ds & = & \int_{\partial M_{j_0}}H_0|P^\mp\Phi|^2ds.
\end{eqnarray*}

Using this relation in (\ref{masse1}) and since $P^\pm\Psi=P^\pm\Phi$ on $\partial M_{j_0}$, we get:
\begin{eqnarray*}
\int_{\partial M_{j_0}}\<{\bf{D}}\Psi,\Psi\>ds\leq\frac{n-1}{2}\int_{\partial M_{j_0}}H_0|\Psi|^2ds
\end{eqnarray*}

with equality if and only if $P^\mp\Psi=P^\mp\Phi$ on $\partial M_{j_0}$. Since we assumed that $0\leq H_0\leq H^{(j_0)}$:
\begin{eqnarray}
\int_{\partial M_{j_0}} \Big(\<{\bf{D}}\Psi,\Psi\>-\frac{n-1}{2}H^{(j_0)}|\Psi|^2\Big)ds\leq 0.
\end{eqnarray}

Now if we look at the boundary term in (\ref{reil}) for $j\neq j_0$, we have:
\begin{eqnarray*}
\int_{\partial M_j} \Big(\<{\bf{D}}\Psi,\Psi\> -\frac{n-1}{2}H^{(j)}|\Psi|^2\Big)ds & = &  -\frac{n-1}{2}\int_{\partial M_j}H^{(j)}|P^\mp\Psi|^2ds
\end{eqnarray*}

because (\ref{prolonge}) gives $P^\pm\Psi=P^\pm\widetilde{\Phi}_j=0$ and so:
\begin{eqnarray}\label{bordj}
\sum_{j\neq j_0}\int_{\partial M_j} \Big(\<{\bf{D}}\Psi,\Psi\>-\frac{n-1}{2}H^{(j)}|\Psi|^2\Big)ds\leq 0
\end{eqnarray}

since $H^{(j)}\geq 0$. Moreover, equality occurs in (\ref{bordj}) if and only if $P^\mp\Psi=0$ (because $H^{(j)}$ is a non zero smooth function on $\partial M_j$). Using (\ref{masse2}) and (\ref{bordj}), we conclude that the boundary term in (\ref{reil}) is nonpositive and so we have equality in the spinorial Reilly formula. Finally we have shown that the spinor field $\Psi\in\Gamma(\Sigma M)$ satisfies:
\begin{eqnarray}\label{concl}
\nabla\Psi=0\quad\text{and}\quad\Psi_{|\partial M}=\widetilde{\Phi}_{|\partial M}.
\end{eqnarray}

In this case, the boundary has to be connected. Indeed, since the spinor $\Psi$ is parallel, it has a non zero constant norm on $M$ (since $M$ is connected), hence on every connected component of $\partial M$. However since $\widetilde{\Phi}_{j}=0$ for $j\neq j_0$ and $\Psi_{|\partial M_{j_0}}=\Phi$, (\ref{concl}) holds only if the boundary is connected (otherwise the norm of $\Psi$ is not constant). On the other hand, using the spinorial Gauss formula (\ref{gauss}) and since $\Psi$ is parallel, we easily check that ${\bf{D}}\Psi=\frac{n-1}{2}H^{(j_0)}\Psi$ and $H^{(j_0)}=H_0$ since $\Psi$ satisfies (\ref{diracequation}) and has no zeros.
\hfill $\square$\\

\begin{remark}\label{impulsion}
Under the assumptions of Theorem \ref{main}, the second fundamental form of $(\partial M,g)$ is completely determined by the spinor field $\Phi\in\Gamma(\mathbf{S})$ which satisfies (\ref{diracequation}) and more precisely by its energy-momentum tensor $T_\Phi$. Indeed, we proved that the spinor $\Phi$ is a generalized Killing spinor field (in the sense of \cite{bgm}), that is it satisfies: 
\begin{eqnarray*}
\nabla^{{\bf S}}_X\Phi=-\frac{1}{2}\gamma^{{\bf S}}\big(A(X)\big)\Phi,
\end{eqnarray*}

for all $X\in\Gamma\big(T(\partial M)\big)$. Thus we can easily check that (see \cite{Mo}):
\begin{eqnarray*}
A(X,Y):=g\big(A(X),Y\big)=2 T_\Phi(X,Y),
\end{eqnarray*}

where $T_\Phi$ is the energy-momentum tensor associated with $\Phi$ defined (on the complement set of zeros of $\Phi$) by:
\begin{eqnarray*}
T_\Phi(X,Y):=\frac{1}{2}{\rm Re}\,\<\gamma(X)\nabla_Y^{{\bf S}}\Phi+\gamma(Y)\nabla_X^{{\bf S}}\Phi,\frac{\Phi}{|\Phi|^2}\>
\end{eqnarray*}

for $X\in\Gamma\big(T(\partial M)\big)$. 
\end{remark}

\begin{remark} 
In the $3$-dimensional case, one can refine the conclusion of Theorem \ref{main}. In fact, if $(M^3, g)$ is a manifold satisfying these assumptions then its Ricci tensor vanishes and we can conclude that $(M^3, g)$ is flat.
\end{remark}

Thanks to Theorem~\ref{main}, we obtain new rigidity results for compact manifolds with boundary which highlight that the boundary behaviour of the metric has an influence on the metric in the interior of the manifold. The main argument is to observe that the Dirac equation (\ref{diracequation}) has a nice geometric interpretation. Indeed, it is quite easy to show that if $(\Sigma^{n-1},g)$ is a smooth oriented hypersurface with mean curvature $H_0$ in a Riemannian spin manifold $(N^n,\widetilde{g})$ carrying a parallel spinor field $\Phi\in\Gamma(\Sigma N)$ then $\Phi_{|\Sigma}$ satisfies the Dirac equation (\ref{diracequation}). As a consequence of this remark, we get a counterpart of results of Ros \cite{ros} and Hang-Wang \cite{hangwang1} in the spinorial setting:
\begin{theorem}\label{geometric}
Let $(M^n, g_1)$ be an $n$-dimensional complete Riemannian spin manifold with nonnegative scalar curvature and $(\Sigma^{n-1}, g)$ a compact hypersurface endowed with the induced Riemannian and spin structures. Assume that there exists an isometric immersion
\begin{eqnarray*}
 F_1:(\Sigma^{n-1}, g)\rightarrow(M^n, g_1)
\end{eqnarray*}

with mean curvature $H_1$ and such that $F_1(\Sigma^{n-1})$ bounds a compact domain $\Omega$ in $M$. Then if there is another isometric immersion
\begin{eqnarray*}
F_2:(\Sigma^{n-1},g)\rightarrow(N^n,g_2)
\end{eqnarray*}

where $(N^n, g_2)$ is a Riemannian manifold carrying a parallel spinor and such that the mean curvature $H_2\geq 0$ of $F_2$ satisfies $H_1\geq H_2$, then the domain $(\Omega,g_1)$ carries a parallel spinor.
\end{theorem}

{\it Proof:} The restriction of the parallel spinor field by $F_2$ yields a solution of the Dirac equation (\ref{diracequation}) with $H_0=H_2$ and thus Theorem~\ref{main} allows to conclude.
\hfill$\square$

\begin{remark}
The simply connected manifolds carrying parallel spinor fields are classified in \cite{mkw} and thus Theorem~\ref{geometric} can be applied for manifolds with boundary whose boundary can be isometrically immersed in this class of manifolds.
\end{remark}

\begin{remark}\label{remeuclidien}
In the statement of Theorem~\ref{geometric}, if we assume that $F_2$ is an isometric immersion in $\mathbb{R}^n$ endowed with its Euclidean metric, one can check that $(\Omega,g_1)$ is flat. If moreover $F_2$ is an isometric embedding, an argument similar to {\rm Proposition} $2$ of \cite{hangwang1} allows to conclude that $(\Omega,g_1)$ is isometric to a domain in $\mathbb{R}^n$.
\end{remark}

As a corollary of Theorems \ref{main} and \ref{geometric}, we give a proof of a conjecture by Schroeder and Strake \cite{schroeder}. More precisely, they prove:
\begin{thrmst}\label{st}
Let $(M^n,g)$ a compact and connected Riemannian manifold with nonnegative Ricci curvature and with convex boundary (that is $A\geq 0$). Assume that the sectional curvature of $M$ vanishes on a neighbourhood of $\partial M$ and that one of the following conditions holds:
\begin{enumerate}
\item $\partial M$ is simply connected
\item the dimension of $\partial M$ is even and $\partial M$ is strictly convex at some point $p\in\partial M$.
\end{enumerate}

Then $(M^n,g)$ is flat.
\end{thrmst}

As pointed out in \cite{schroeder}, the condition on the sectional curvature is very strong and the authors conjecture that their results should hold under the weaker condition of vanishing of the sectional curvature {\it along} the boundary. In \cite{hangwang1}, F. Hang and X. Wang proved the part $(1)$ of this conjecture. In fact, they observe that it is enough to impose the nonnegativity of the mean curvature of the boundary and not necessarily its convexity. If the manifold is spin, they can relax the condition on the Ricci curvature by only assuming the nonnegativity of the scalar curvature. However, in this case they need the convexity of the boundary. Here the spin assumption is essentially technical because their proof relies on some positive mass theorems (see \cite{shitam}) proved with spinors \cite{witten}. We give here a generalization of Hang and Wang's result and thus of a part of the Schroeder and Strake's conjecture. More precisely, we get:
\begin{corollary}\label{schroeder}
Let $(M^n, g)$ an $n$-dimensional compact, connected Riemannian spin manifold with boundary and with nonnegative scalar curvature. If every component of the boundary of $M$ is simply connected with nonnegative mean curvature and the sectional curvature vanishes on $\partial M$, then the boundary has only one connected component and $(M^n, g)$ is flat.
\end{corollary}

{\it Proof:} We first remark that since the sectional curvature $\kappa^M$ of $M$ is identically zero on $\partial M$, the Weingarten map $A$ satisfies the Gauss and Codazzi equations:
\begin{eqnarray*}
(\nabla^{\partial M}_X A)Y & = & (\nabla^{\partial M}_Y A)X\\
R^{\partial M}(X,Y)Z & = & g\big(A(Y),Z\big)A(X)-g\big(A(X),Z\big)A(Y),
\end{eqnarray*}

for $X,Y,Z\in\Gamma\big(T(\partial M)\big)$. On the other hand, the boundary $\partial M$ is simply connected, then the fundamental theorem for hypersurfaces (see \cite{KN69} for example) ensures the existence of an isometric immersion $F$ of $(\partial M, g)$ in $(\mathbb{R}^n,eucl)$ with Weingarten map given by $A$. With this immersion, we get $2^{[n/2]}$ spinor fields $\Phi_i\in\Gamma ({\bf{S}})$ such that:
\begin{eqnarray*}
\nabla^{{\bf{S}}}_X\Phi_i=-\frac{1}{2}\gamma^{\bf{S}}\big(A(X)\big)\Phi_i
\end{eqnarray*}

for all $X\in\Gamma\big(T(\partial M)\big)$ and thus ${\bf{D}}\Phi_i=\frac{n-1}{2} H\Phi_i$. These spinor fields are the restriction (by $F$) on ${\bf S}$ of $2^{[n/2]}$ parallel spinor fields on $\Sigma\mathbb{R}^n$. Thus the assumptions of Theorem~\ref{main} (or Theorem~\ref{geometric}) are fulfilled and one concludes
that the boundary is connected and that each spinor field $\Phi_i$ comes from a parallel spinor field on $M$. Finally, we get a maximal number of parallel spinor fields and thus $(M^n,g)$ is flat.
\hfill $\square$\\

Another application of Theorems \ref{main} and \ref{geometric} is given by a simple proof of a rigidity result for the unit Euclidean ball. This result has been proved by P. Miao \cite{miao} as a consequence of a positive mass theorem for asymptotically flat manifolds for which the metric is not smooth along a hypersurface. The proof we give here relies only on Spin Geometry and does not use this strong but quite technical result. We show:
\begin{corollary}\label{boule1}
Let $(M^n, g)$ be a compact and connected Riemannian spin manifold with smooth boundary. Assume that the scalar curvature of $M$ is nonnegative and that the boundary is isometric to the standard sphere $\mathbb{S}^{n-1}$ with mean curvature satisfying $H\geq 1$. Then $(M^n,g)$ is isometric to the unit ball of $\mathbb{R}^n$.
\end{corollary}

{\it Proof:} From Theorem~\ref{main} (or \ref{geometric}), we have a basis of $\Sigma M$ made of parallel spinor fields whose restrictions correspond to extrinsic Killing spinors on the boundary (in the sense of \cite{HM}). On the other hand, using the spinorial Gauss formula (\ref{gauss}) we see that the boundary has to be totally umbilical with constant mean curvature. As a conclusion, $M$ is a compact flat Riemannian spin manifold whose boundary is a totally geodesic round sphere, then $(M^n, g)$ is isometric to the unit Euclidean ball of $(\mathbb{R}^n, eucl)$.
\hfill $\square$\\

\begin{remark}
One can check that the assumptions of Corollary~\ref{boule1} are not covered in the work of Hijazi and Montiel \cite{HM}. 
\end{remark}

\begin{remark}
It is clear that Corollary~\ref{boule1} holds if the boundary is isometric to the sphere $\mathbb{S}^{n-1}(r)$ with radius $r>0$ and with mean curvature satisfying $H\geq 1/r$. In this case, $(M^n,g)$ is isometric to the Euclidean ball with radius $r$. 
\end{remark}

Under the assumptions (and notations) of Theorem \ref{geometric}, one can also ask the following question: can the domain $(\Omega,g_1)$ be isometrically immersed in $(N^n,g_2)$? We don't give an answer to this question here but give some ideas for further investigations. For this, we assume that the image of $(\Sigma^{n-1},g)$ by $F_2$ also bounds a compact domain in $(N^n,g_2)$. With the help of Remark \ref{impulsion}, one concludes that if $A_i$ denotes the second fundamental form of $\Sigma$ for $F_i$ with $i=1,2$, we get $A_1=A_2$. Using the recent work of M.T. Anderson and M. Herzlich \cite{andherz} on the unique continuation properties for Einstein manifolds with boundary, we get:
\begin{corollary}
Let $(\Omega^n,g_1)$ a compact and connected Riemannian spin manifold with nonnegative scalar curvature. Assume that its boundary $\Sigma^{n-1}$ with mean curvature $H_1$ can be isometrically embedded in a Riemannian spin manifold $(N^n,g_2)$ carrying a parallel spinor field with mean curvature less than $H_1$. Then there exists a neighborhood of $(\Sigma^{n-1},g)$ in $(\Omega^n,g_1)$ which can be isometrically embedded in $(N^n,g_2)$.
\end{corollary}


\section{Domains with negative scalar curvature}\label{hyperbolique}


In this section, we prove rigidity results similar to Theorems \ref{main} and \ref{geometric} in the hyperbolic setting. More precisely, we consider an $n$-dimensional connected and compact Riemannian spin manifold  $(M^n, g)$  with smooth boundary $\partial M$. We also assume that the scalar curvature of $M$ (for the metric $g$) satisfies $R\geq -n(n-1)$ and that the mean curvature $H$ of the boundary $\partial M$ is nonnegative (and non identically zero). \\

Here it is important to note that the proof of Theorem~\ref{main} lies on two important facts: the Schr\"odinger-Lichnerowicz formula and a suitable boundary condition for the Dirac operator of $M$. So we first recall the hyperbolic version of the Schr\"odinger-Lichnerowicz formula where a proof can be found in \cite{anderson}, \cite{HMR} or \cite{minoo}:
\begin{eqnarray}\label{reillyhyp}
\int_{M}\big(|P\psi|^2+\frac{1}{4}\widetilde{R}|\psi|^2-\frac{n-1}{n}|\widetilde{D}^\pm\psi|^2\big) dv 
=\int_{\partial M} \big(\<\widetilde{\bf{D}}^\pm\psi,\psi\> -\frac{n-1}{2}{H}|\psi|^2\big)ds
\end{eqnarray}

where $\widetilde{R}:=R+n(n-1)$, $\widetilde{D}^\pm:=D\mp\frac{n}{2}i$ and $\widetilde{\bf{D}}^\pm:={\bf{D}}\pm\frac{n-1}{2}i\gamma(\nu)$. The operator $P$ in (\ref{reillyhyp}) is the twistor operator (or Penrose operator) locally given by:
\begin{eqnarray}\label{twisteur}
P_X\psi:=\nabla_X\psi+\frac{1}{n}\gamma(X)D\psi,
\end{eqnarray}

for $\psi\in\Gamma(\Sigma M)$ and $X\in\Gamma(T M)$. A spinor $\psi$ such $P\psi=0$ is called a twistor spinor. On the other hand, one can check that the operator $\widetilde{\bf{D}}^\pm$ which appears in the boundary term of (\ref{reillyhyp}) is an elliptic first order self-adjoint differential operator and its spectrum is an unbounded sequence of real numbers.  Observe now that the choice of the boundary condition deeply lies on its behaviour with respect to the {\it twisted} Dirac operator $\widetilde{\bf{D}}^\pm$. The condition used in Section~\ref{euclidien} is not appropriate and that is why we will use another elliptic boundary condition for the Dirac operator $D$: {\it the condition associated with a chirality operator}.
This kind of condition does not exist on all manifolds since it needs a chirality operator, that is a linear map: 
\begin{eqnarray*}
 G:\Gamma(\Sigma M)\longrightarrow\Gamma(\Sigma M),
\end{eqnarray*}

such that:
\begin{align}\label{chi}
G^2 &=Id,\quad\<G\psi,G\varphi\>=\<\psi,\varphi\>\\
\nabla_{X}(G\psi) &=G(\nabla_{X}\psi),\quad\gamma(X)G(\psi)=-G(\gamma(X)\psi),
\end{align}

\noindent for all $X\in\Gamma(T M)$ and for all spinor fields $\psi,\varphi\in\Gamma(\Sigma M)$. If we assume the existence of such an operator, we can define an involution on ${\bf{S}}$ by:
\begin{eqnarray*}
\gamma(\nu)G:\Gamma({\bf S})\longrightarrow\Gamma({\bf S}),
\end{eqnarray*}

which gives a decomposition of the spinor bundle ${\bf{S}}$ into the direct sum of two eigensubbundle associated with the eigenvalues $1$ and $-1$. The pointwise orthogonal projection:
\begin{eqnarray*}
B^\pm:=\frac{1}{2}(Id\pm\gamma(\nu)G)
\end{eqnarray*}

on the eigensubbundle associated with the eigenvalue $\pm 1$ defines an elliptic boundary condition for the Dirac operator $D$ of $M$. For more details on this boundary condition, we refer to \cite{hijazi.montiel.roldan:01} or \cite{mathese} for example.\\

We can now state the main result of this section which can be seen as a hyperbolic counterpart of Theorem~\ref{main} of Section~\ref{euclidien}. Indeed, we have:
\begin{theorem}\label{main1}
Let $(M^n, g)$ a connected and compact Riemannian spin manifold with smooth boundary equipped with a chirality operator $G$. Assume that the scalar curvature of $M$ satisfies $R\geq -n(n-1)$ and that every connected component $\partial M_j$ of $\partial M$ in $M$ has nonnegative mean curvature $H^{(j)}\geq 0$. If there exists a smooth spinor field $\Phi\in\Gamma ({\bf{S}}_{j_0})$ such that:
\begin{eqnarray}\label{diracequationhyp}
{\widetilde{\bf{D}}}^\pm\Phi=\frac{n-1}{2} H_0\Phi,
\end{eqnarray}

where $H_0$ is a nonnegative (and non zero) smooth function on $\partial M$ with $0\leq H_0\leq H^{(j_0)}$, the manifold $(M^n, g)$ carries an imaginary Killing spinor, the boundary is connected and $H^{(j_0)}=H_0$.
\end{theorem} 

In order to prove this result, we first show the following lemma:
\begin{lemma}\label{lemchi}
Under the assumptions of Theorem~\ref{main1}, the Dirac operator with domain:
\begin{eqnarray*}
\widetilde{D}^\pm:\{\psi\in H^2_1(\Sigma M)\,:\,
B^\pm\psi_{|\partial M}=0\}\longrightarrow L^2(\Sigma M)
\end{eqnarray*}

\noindent is invertible.
\end{lemma}

{\it Proof:} Suppose that there exists a non trivial spinor field $\varphi_0\in\Gamma(\Sigma M)$ solution of the boundary value problem:
$$\left\lbrace
\begin{array}{ll}
\widetilde{D}^+\varphi_0=0 & \qquad\text{on}\;M\\
B^\pm\varphi_{0\,|\partial M}=0 & \qquad\text{along}\;\partial M
\end{array}
\right.$$

that is:
$$\left\lbrace
\begin{array}{ll}
D\varphi_0=\frac{n}{2}i\varphi_0 & \qquad\text{on}\;M\\
B^\pm\varphi_{0\,|\partial M}=0 & \qquad\text{along}\;\partial M.
\end{array}
\right.$$

The Green formula gives for all $\psi\in\Gamma(\Sigma M)$:
\begin{eqnarray}\label{ipp}
\int_M\<D\psi,\psi\>dv-\int_M\<\psi,D\psi\>dv=-\int_{\partial M}\<\gamma(\nu)\psi,\psi\>ds
\end{eqnarray}

and by sesquilinearity of the Hermitian product on $\Sigma M$, we get for $\psi=\varphi_0$:
\begin{eqnarray*}
\int_M|\varphi_0|^2dv=0,
\end{eqnarray*}

which implies that $\varphi_0\equiv 0$ on $M$ and so a contradiction with our assumptions. The conclusion follows from the fact that:
\begin{eqnarray*}
(\widetilde{D}^+)^*=\widetilde{D}^-\,\,\Longrightarrow\,\,{\rm CoKer}\,(\widetilde{D}^+)\simeq {\rm Ker}\,(\widetilde{D}^-)=\{0\},
\end{eqnarray*}

where we used (\ref{ipp}).
\hfill$\square$\\

In the following lemma, we study the behaviour of the twisted Dirac operator $\widetilde{\bf{D}}^+$ with respect to the boundary condition associated with a chirality operator.
\begin{lemma}\label{lemchi1}
If $\Phi\in\Gamma ({\bf S})$ is a smooth spinor field satisfying the Dirac equation (\ref{diracequationhyp}), we have:
\begin{enumerate}
\item $\widetilde{\bf{D}}^+(B^\pm\Phi)=\frac{n-1}{2}H_0B^\mp\Phi$
\item $\int_{\partial M}H_0 |B^\pm\Phi|^2ds=\int_{\partial M}H_0 |B^\mp\Phi|^2ds$
\end{enumerate}
\end{lemma}

{\it Proof:}
For $(1)$, it is sufficient to note that:
\begin{eqnarray*}
\widetilde{\bf{D}}^+(B^\pm\psi)=B^\mp(\widetilde{\bf{D}}^+\psi),
\end{eqnarray*}

and since $\Phi$ is a solution of (\ref{diracequationhyp}), a simple identification of the components of a spinor field with respect to the decomposition associated with the orthogonal projections $B^\pm$ gives the result. Point $(2)$ follows from $(1)$ and from the symmetry of the operator $\widetilde{\bf{D}}^+$.
\hfill$\square$\\

The proof of Theorem~\ref{main1} is then similar to the one of Theorem~\ref{main}. We don't give the details here since Lemmas~\ref{lemchi} and \ref{lemchi1} show that the condition associated with a chirality operator has the same behaviour with respect to the twisted Dirac operator $\widetilde{\bf{D}}^\pm$ as the ${\rm MIT}$ boundary condition with respect to the extrinsic Dirac operator ${\bf{D}}$.\\

As a consequence of Theorem~\ref{main1}, we get rigidity results for manifolds with boundary with scalar curvature bounded by below by a negative constant. The main geometric result we prove in this setting is:
\begin{theorem}\label{hypgeometric}
Let $(M^n, g_1)$ be an $n$-dimensional complete Riemannian spin manifold with scalar curvature satisfying $R\geq -n(n-1)$ and let $(\Sigma^{n-1},g)$ be a compact hypersurface endowed with the induced Riemannian and spin structures. Assume that there exists an isometric immersion
\begin{eqnarray*}
 F_1:(\Sigma^{n-1},g)\rightarrow(M^n,g_1)
\end{eqnarray*}

with mean curvature $H_1$ and such that $F_1(\Sigma^{n-1})$ bounds a compact domain $\Omega$ in $M$. Then if there is another isometric immersion
\begin{eqnarray*}
F_2:(\Sigma^{n-1},g)\rightarrow(N^n,g_2)
\end{eqnarray*}

where $(N^n,g_2)$ is a Riemannian manifold carrying an imaginary Killing spinor (with constant $\mp (i/2)$) and such that the mean curvature $H_2\geq 0$ of $F_2$ satisfies $H_1\geq H_2$, then the domain $(\Omega,g_1)$ carries an imaginary Killing spinor (with same constant).
\end{theorem}

{\it Proof:} It is sufficient to note that the immersion $F_2$ of $(\Sigma^{n-1},g)$ in $(N^n,g_2)$ yields to the existence of a smooth spinor field on ${\bf S}$ satisfying the Dirac equation (\ref{diracequationhyp}). Moreover the assumptions on the mean curvatures $H_1$ and $H_2$ enable us to apply Theorem \ref{main1} and thus one concludes that the domain $(\Omega,g_1)$ carries an imaginary Killing spinor.
\hfill$\square$\\

From this result, we obtain a hyperbolic version of Corollary \ref{schroeder} which gives a counterpart of Schroeder and Strake's conjecture in this setting. More precisely, we have: 
\begin{corollary}\label{schroederhyp}
Let $(M^n, g)$ be a connected and compact Riemannian spin manifold such that the scalar curvature satisfies $R\geq -n(n-1)$. If every component of the boundary of $M$ is simply connected with nonnegative mean curvature and the sectional curvature is $-1$ on $\partial M$, then the boundary has only one connected component and $(M^n, g)$ is hyperbolic.
\end{corollary}

{\it Proof:} Since $\partial M$ is simply connected and the sectional curvature $\kappa^M$ of $M$ is $-1$ on $\partial M$, the boundary $(\partial M, g)$ can be isometrically immersed in the standard hyperbolic space $(\mathbb{H}^n,g_{st})$ which is endowed with a maximal number of imaginary Killing spinors. Thus using Theorem \ref{hypgeometric} (or Theorem \ref{main1}), we get the existence of a maximal number of imaginary Killing spinors on $(M^n,g)$ and thus by \cite{h1} and \cite{h2} it has to be hyperbolic.
\hfill$\square$\\

\begin{remark}
All the preceding results hold for even dimensional manifolds since in these dimensions, the spinor bundle is endowed with a chirality operator (the volume element of the spinor bundle).
\end{remark}

\begin{remark}
It is clear that with Theorem \ref{main1}, we can prove a rigidity result for geodesic balls of the hyperbolic space $\mathbb{H}^n$. However in this case we need the existence of a chirality operator which is in fact not necessary as explained below. Indeed, one can obtain such a result from an estimate on the first eigenvalue of the twisted Dirac operator $\widetilde{\bf{D}}^\pm$ proved by Hijazi, Montiel and Rold\'an in \cite{HMR}. In fact, they show that if $(M^n, g)$ is a connected compact Riemannian spin manifold with smooth boundary $\partial M$ such that the scalar curvature of $M$ is bounded from below by $-n(n-1)$ and the mean curvature $H$ of $\partial M$ in $M$ is nonnegative:
\begin{eqnarray}\label{hijmonrol}
\lambda_1^\pm\geq \frac{n-1}{2}\inf_{\partial M}H
\end{eqnarray}

\noindent where $\lambda_1^\pm$ denotes the first eigenvalue of $\widetilde{\bf{D}}^\pm$. Moreover equality occurs if and only if the eigenspinors associated with the eigenvalue $\lambda_1^\pm$ consist of restrictions to $\partial M$ of imaginary Killing spinors on $M$. With this result, we can prove:\\

{\rm ``Let $(M^n,g)$ a compact domain in a complete Riemannian spin manifold with scalar curvature satisfying $R\geq -n(n-1)$. Assume that the boundary of $M$  is isometric to the standard round sphere $\mathbb{S}^{n-1}\big(\frac{1}{\sqrt{\alpha^2-1}}\big)$ ($\alpha>1$) and that its mean curvature is such that $H\geq \alpha$. Then $(M^n, g)$ is isometric to the standard ball of $\mathbb{H}^n$ whose boundary is totally umbilical (and isometric to the standard sphere of radius $\frac{1}{\sqrt{\alpha^2-1}}$).''}\\

\noindent Indeed, it is enough to note that for each real Killing spinor on the boundary (which exists since the boundary is isometric to a round sphere), we can construct an eigenspinor for the Dirac operator $\widetilde{\bf{D}}^\pm$ associated with the eigenvalue $\frac{n-1}{2}\alpha$. More precisely, if $\Psi$ denotes a real Killing spinor with Killing number $\frac{1}{2\sqrt{\alpha^2-1}}$ then the spinor field defined by:
\begin{eqnarray*}
\Psi^\pm:=\Psi\pm(\alpha-\sqrt{\alpha^2-1})i\gamma(\nu)\Psi,
\end{eqnarray*}

\noindent satisfies:
\begin{eqnarray}\label{vp}
\widetilde{{\bf{D}}}^\pm\Psi^\pm=\frac{n-1}{2}\alpha\Psi^\pm.
\end{eqnarray}

\noindent Thus the assumption $H\geq \alpha$ shows that the equality case in (\ref{hijmonrol}) is reached and then from \cite{HMR} the spinor field $\Psi^\pm$ is the restriction of an imaginary Killing spinor on $(M^n,g)$. Since there is a maximal number of real Killing spinors on $\partial M$, one easily check that there is a maximal number of imaginary Killing spinors on $M$ and \cite{h1} and \cite{h2} give the result.
\end{remark}



\end{document}